\documentclass[a4paper,reqno,11pt]{amsart}

\usepackage[all]{xy}
\usepackage[english]{babel}
\usepackage{amssymb}
\usepackage{graphicx}
\usepackage{hhline}
\usepackage{xcolor}
\usepackage{shuffle}
\usepackage{wasysym}
\usepackage{amsmath,amsthm,amssymb,xypic,verbatim}

\usepackage{geometry}
\newcommand{\p}{\mathbb{P}}

\newcommand{\N}{\mathbb{N}}

\newcommand{\C}{\mathbb{C}}
\newcommand{\R}{\mathbb{R}}
\newcommand{\X}{\mathbb{X}}

\newcommand{\LL}{\mathcal{L}}
\newcommand{\G}{\mathcal{G}}
\newcommand{\Lie}{\mathop{\rm Lie}\nolimits} 
\newcommand{\f}{\varphi}
\newcommand{\V}{\mathcal{V}}
\newcommand{\RR}{\mathcal{R}}
\newcommand{\AAA}{\mathcal{A}}
\newcommand{\U}{\mathcal{U}}

\newcommand{\Sym}{\mathop{{\rm Sym}}\nolimits}

\newcommand{\sgn}{\mathop{\rm sgn}\nolimits}
\newcommand{\Pf}{\mathop{\rm Pf}\nolimits}

\theoremstyle{plain}

\newtheorem{thm}{Theorem}
\newtheorem{pro}[thm]{Proposition}
\newtheorem{lem}[thm]{Lemma}   
\newtheorem{cor}[thm]{Corollary}

\theoremstyle{definition}

\newtheorem{notaz}[thm]{Notation}

\newtheorem{es}[thm]{Example}

\newtheorem{dfn}[thm]{Definition}

\newtheorem{rmk}[thm]{Remark}

\title{The Rough Veronese variety}
\author{Francesco Galuppi}
\address{Max Planck Institute for Mathematics in the Sciences\\
	Inselstra\ss{}e 22\\
	04103 Leipzig\\
	Germany}
\email{galuppi@mis.mpg.de}
\date{}
\subjclass{14Q15, 14M25, 60H99}
\keywords{Rough Veronese variety, signature tensor, rough path, Lie algebra, Lyndon words, exponential}
\thanks{The Nonlinear Algebra research group at Max Planck Institute is a very stimulating environment, and I would like to thank Rosa Prei\ss{} for introducing me to the topic, Jane Coons for her help with the combinatorics, Joscha Diehl for his many explanations about the tensor algebra, Mateusz Micha\l{}ek for his patience and his insight, Laura Colmenarejo for providing lots of examples, and Max Pfeffer, Anna Seigal and Bernd Sturmfels for useful suggestions and careful proofreading.}

\begin{document}
	\maketitle

\begin{abstract}
	We study signature tensors of paths from an algebraic geometric viewpoint. The signatures of a given class of paths parametrize a variety inside the space of tensors, and these signature varieties provide both new tools to investigate paths and new challenging questions about their behavior.
	This paper focuses on signatures of rough paths. Their signature variety shows surprising analogies with the Veronese variety, and our aim is to prove that this so-called Rough Veronese is toric. The same holds for the universal variety.
	Answering a question of Am\'{e}ndola, Friz and Sturmfels, we show that the ideal of the universal variety does not need to be generated by quadrics.
\end{abstract}
	
\section*{Introduction}\label{sec:intro}

A \textit{path} is a map $X:[0,1]\to\R^d$. Classically, the components $X_1,\dots,X_d$ of $X$ are assumed to be sufficiently smooth. 
 For every positive integer $k$, it is therefore possible to define an order $k$ tensor $\sigma^{(k)}(X)\in(\R^d)^{\otimes k}$, whose $(i_1\ldots i_k)$-th entry is
\[\int_{0}^{1}\int_{0}^{t_k}\dots\int_{0}^{t_3}\int_{0}^{t_2}\dot{X}_{i_1}(t_1)\cdot\ldots\cdot\dot{X}_{i_k}(t_k)dt_1\dots dt_k.
\]
By convention, we define $\sigma^{(0)}(X)=1$. The sequence
\[\sigma(X)=(\sigma^{(k)}(X)\mid k\in\N)\]
is called the \textit{signature} of $X$. Signatures were first defined in \cite{Chen}, and they enjoy many interesting properties. For instance, the signature allows to uniquely recover a sufficiently smooth path up to a mild equivalence relation (see \cite[Theorem 4.1]{chenuniqueness}).

Many physical behaviors and experiments can be modeled by using paths, and signatures are useful tools to encode in a compact form the information carried by paths. In \cite{AFS}, the authors consider signature tensors from an algebraic geometry perspective. If we fix a certain class of paths and the order $k$ of the tensors, then the $k$-th signature $\sigma^{(k)}$ is a polynomial map into $(\R^d)^{\otimes k}$. Its image variety parametrizes the closure of the set of all $k$-th signatures of paths of the chosen class. The study of this map, this \textit{signature variety} and its geometric properties is interesting for many reasons.
For instance, in applied problems one sometimes has a signature, coming from empirical data, and wants to know if there is a path of a certain kind (say, piecewise linear) having that given signature. Knowing whether the map is injective, or at least finitely many to one, tells us if there are chances to solve this inverse problem. Another issue is the study of singularities. What does it mean, for a path, to have a signature which is a singular point in the variety?

In \cite{AFS} we find a detailed study of the signature varieties of polynomial paths, piecewise linear paths and also random paths arising from Brownian motion. These three classes of paths have a common generalization: the class of rough paths. Rough paths %, introduced in \cite{origin of rough paths}, 
have a number of applications, for instance the study of controlled ODEs and stochastic PDEs (see \cite{Friz}), as well as sound compression (see \cite{soundcompression}).
While they are not necessarily smooth, it is possible to define their signature, and therefore to study their signature variety. Even at a first glance, such a variety exhibits analogies with the Veronese variety, and it is therefore named the \textit{Rough Veronese variety} in \cite[Section 5.4]{AFS}. The main purpose of this paper is to study its geometry. We will prove that the Rough Veronese variety is toric, and we will characterize the monomials parameterizing it.

Sometimes it is useful to perform explicit computations in order to handle concrete examples. In these cases, we use the software Macaulay2, freely available at {\tt{www.math.uiuc.edu/Macaulay2}}.

\section*{Notations and preliminaries}\label{sec:notation}
The $k$-th signature of a path $X$ belongs to $(\R^d)^{\otimes k}$, but we need a space to store the whole signature $\sigma(X)$. In this Section we define such a space, which has a rich algebraic structure. We will recall the features we need, but we do not attempt to any extent to describe all its properties. Every definition and result of this section can be found in \cite{bible}.

\begin{dfn}
	The \textit{tensor algebra} over $\R^d$ is the graded $\R$-vector space
	\[T((\R^d))= %\bigoplus_{k\in\N}(\R^d)^{\otimes k}	=
	\R\times\R^d\times(\R^d)^{\otimes 2}\times\ldots\]
	of formal power series in the non-commuting variables $x_1,\dots,x_d$. It is an $\R$-algebra with respect to the tensor product, and we denote by $p_k:T((\R^d))\to(\R^d)^{\otimes k}$ the projection. The algebraic dual of $T((\R^d))$ is the graded algebra
	\[T(\R^d)=\R\langle x_1,\dots,x_d\rangle\]
	of polynomials in the non-commuting variables $x_1,\dots,x_d$. It is the free $\R$-algebra over $x_1,\dots,x_d$.
\end{dfn}

\begin{notaz} Given an element $T\in T((\R^d))$, we denote by $T_{i_1\dots i_k}$ the $(i_1\dots i_k)$-th entry of the order $k$ element of $T$. For $y\in\R$, we denote $T_y((\R^d))=\{T\in T((\R^d))\mid T_0=y\}$. 
	
It will often be convenient to identify a degree $k$ monomial $x_{i_1}\cdot\ldots\cdot x_{i_k}$ with the word $w={\bf i_1}\ldots {\bf i_k}$ in the alphabet $\{{\bf 1},\dots,{\bf d}\}$. The number $k$ is called the \textit{length} of $w$ and it is denoted by $|w|$. The degree 0 monomial corresponds to the empty word ${\bf e}$. We will write letters in bold in order to distinguish the number 1 from the letter ${\bf 1}=x_1$. In this way, the product of two words $v$ and $w$ is simply the word obtained by writing $v$ followed by $w$. It is called the \textit{concatenation product} and sometimes denoted by $v\cdot w$. The natural duality pairing
	\[\langle-,-\rangle:T((\R^d))\times T(\R^d)\to\R
	\]
is given by $\langle T,{\bf i_1}\ldots {\bf i_k}\rangle=T_{i_1\dots i_k}$, extended by linearity.
\end{notaz}

Besides the concatenation of words, there is another product on $T(\R^d)$. It will play a very important role in this paper.

\begin{dfn}
	The \textit{shuffle product} of two words $v$ and $w$, denoted by $v\shuffle w$, is the sum of all order-preserving interleavings of them. More precisely, we define the shuffle recursively. If ${\bf i}$ and ${\bf j}$ are letters, then
	\begin{align*}
	&{\bf e} \shuffle w = w \shuffle {\bf e} = w \text{ and }\\
	&(v\cdot {\bf i}) \shuffle (w\cdot {\bf j}) = \left(v\shuffle (w\cdot {\bf j})\right)\cdot {\bf i} + \left((v\cdot {\bf i} )\shuffle w\right)\cdot {\bf j}.
	\end{align*}
			As an example,
	\[{\bf 12}\shuffle{\bf 34}={\bf 1234}+{\bf 1324}+{\bf 1342}+{\bf 3124}+{\bf 3142}+{\bf 3412}. \]
	Again, the shuffle product is extended by linearity to $T(\R^d)$. We will sometimes use the notation
	\[v^{\shuffle n}=\underbrace{v\shuffle\ldots\shuffle v}_{n\mbox{\begin{Small}
			times\end{Small} }}.\]

\end{dfn}

Despite its apparently complicated definition, the shuffle product enjoys good properties. For instance, the space $(T(\R^d),\shuffle,{\bf e})$ is a commutative algebra. Moreover, shuffle behaves nicely with respect to the signatures. In order to be precise, we need to specify the class of paths we deal with. Intuitively, we want to enlarge the class of smooth paths to comprehend those paths for which the iterated integrals are well-defined.

\begin{dfn}
	A path $X:[0,1]\to\R^d$ is \textit{regular} if $\dot{X}:[0,1]\to\R^d$ is well defined, continuous and non-vanishing on $[0,1]$. It is \textit{piecewise regular} if there are $0=a_0<a_1<\ldots<a_{r-1}<a_r=1$ such that $X_{|[a_{i-1},a_i]}$ is regular on $[a_{i-1},a_i]$ for every $i\in\{1,\dots,r\}$. In particular, a piecewise regular path is piecewise continuous and its finitely many discontinuity points are jumps.
\end{dfn}

The class of piecewise regular path is the one considered by Chen in \cite{chenuniqueness}. These paths have a nice interplay with the shuffle product.

\begin{lem}[Shuffle identity]\label{lem:shuffle identity}
	If $X:[0,1]\to\R^d$ is piecewise regular, then
	\[\langle \sigma(X),v\rangle\cdot\langle \sigma(X),w\rangle=\langle \sigma(X),v\shuffle w\rangle
	\]
	for all words $v,w\in T(\R^d)$.
\end{lem}
The shuffle identity is proved in \cite[Proof of Corollary 3.5]{bible} for a more general class of paths. Hence signatures do not fill the whole tensor algebra $T((\R^d))$, but rather they live in the subset of elements with constant term 1 and satisfying the shuffle identity. This is one of the many possible motivations for the next definition.

\begin{dfn}
	We will denote
	\[\G(\R^d)=\{T\in T_1((\R^d))\mid \langle T,v\rangle\cdot\langle T,w\rangle=\langle T,v\shuffle w\rangle
	\mbox{ for all words } v,w\in T(\R^d)\}.\]
\end{dfn}

The object we have just defined is worth a few remarks. It not only contains the signatures of all piecewise regular paths, but it is also a group with respect to the tensor product. This is why its elements are sometimes called \textit{group-like} elements in the literature. $\G(\R^d)$ is not linear, but it is closely related to a linear space.

\begin{dfn}
	On $T((\R^d))$ there is a Lie bracketing $[T,S]=TS-ST$. Then we can define $\Lie(\R^d)$ to be the free Lie algebra generated by $x_1,\dots,x_d$, that is, the smallest linear subspace of $T((\R^d))$ that contains $x_1,\dots,x_d$ and is closed with respect to the bracketing.
\end{dfn}

This Lie algebra and $\G(\R^d)$ are linked by two maps.

\begin{dfn}
	Define $\exp:T_0((\R^d))\to T_1((\R^d))$ by the formal power series
	\[\exp(T)=\sum_{n=0}^{\infty}\frac{T^{\otimes n}}{n!}.\]
	Not surprisingly, $\exp$ has a two-sided inverse $\log:T_1((\R^d))\to T_0((\R^d))$ defined by
	\[\log(S)=\sum_{n=1}^\infty\frac{(-1)^{n+1}}{n}(S-1)^{\otimes n}.\]
\end{dfn}
For our purposes, we need to point out that all the definitions we have recalled have a truncated version. Namely, one can fix $m\in\N$ and consider
\[T^m(\R^d)=\bigoplus_{k=0}^m(\R^d)^{\otimes k},\]
where tensors of order greater than $m$ are set to zero. Inside $T^m(\R^d)$ there are $\G^m(\R^d)$ and $\Lie^m(\R^d)$. The maps $\exp$ and $\log$ are defined in the same way and they restrict to a bijection between $\Lie^m(\R^d)$ and $\G^m(\R^d)$. In order to avoid confusion, we will write $\exp^{(m)}$ to denote the map $T_0^m((\R^d))\to T_1^m((\R^d))$.

We need a last definition before we move to our Rough Veronese variety. % We have to describe a basis for $\Lie^n(\R^d)$.

\begin{dfn}
	A non-empty word $w$ is a \textit{Lyndon word} if, whenever we write $w=pq$ as the concatenation of two nonempty words, we have $w<q$ in the lexicographic order. We denote by $W_{d,m}$ the set of Lyndon words of length at most $m$ in the alphabet $\{{\bf 1},\dots,{\bf d}\}$. 
	
	For every word $w$, there exists a unique pair $(p,q)$ of nonempty words such that $w=pq$ and $q$ is minimal with respect to lexicographic order. The \textit{bracketing} of $w$ is $[p,q]=pq-qp$.
\end{dfn}

We care about Lyndon words because $\Lie^m(\R^d)$ has a basis consisting of all bracketings of Lyndon words of length at most $m$. Recall that the M\"{o}bius function $\mu:\N\to\N$ sends a natural number $t$ to
\[\mu(t)=\begin{cases}
0 & \mbox{ if $t$ is divisible by the square of a prime},\\
1 & \mbox{ if $t$ is the product of an even number of distinct primes},\\
-1 & \mbox{ if $t$ is the product of an odd number of distinct primes}.\\
\end{cases}
\]
Then the number of length $l$ Lyndon words in the alphabet $\{{\bf 1},\dots,{\bf d}\}$, denoted by $\mu_{l,d}$, is
\[\mu_{l,d}=\sum_{t\mid l}\frac{\mu(t)}{l}d^{\frac{l}{t}}\]
and therefore, as a vector space, $\Lie^m(\R^d)$ has dimension
\[\dim\Lie^m(\R^d)=\sum_{l=1}^m \sum_{t\mid l}\frac{\mu(t)}{l}d^{\frac{l}{t}}.\]

\section*{Signatures of rough paths}\label{sec:enter roughveronese}
In stochastic analysis, people often deal with paths that are far from being smooth, not even piecewise regular. Therefore we want to further expand our class of paths. Given their capital importance in stochastics, it is natural to consider rough paths. In this section we recall the basic definitions and we introduce their signature varieties. The main reference for rough paths is \cite{Multidim}.

Consider a  piecewise regular path $X$ and let $t\in [0,1]$. In the definition of $k$-th signature we can replace the integral on $[0,1]$ with the integral on $[0,t]$. This is the same as restricting $X$ to the sub-interval $[0,t]$, hence we will denote it as $\sigma^{(k)}(X_{|[0,t]})$. As an example,
\[\sigma^{(1)}(X_{|[0,t]})_i=\int_{0}^{t}\dot{X_i}(\lambda)d\lambda=X_i(t)-X_i(0).\]
For every $k$, we notice that $\sigma^{(k)}(X_{|[0,t]})$, as a function of $t$, is a path $[0,1]\to(\R^d)^{\otimes k}$. If we look at the full signature $\sigma(X_{|[0,t]})$, we get a path $[0,1]\to\G(\R^d)$ such that its endpoint is $\sigma(X)$. Moreover, this $\G(\R^d)$-valued path satisfies a H\"{o}lder-like inequality. We will use the symbol $f(t)\apprle g(t)$ to indicate that there is a constant $c$ such that $f(t)\le c\cdot g(t)$ for every $t$.

\begin{lem} Let $X:[0,1]\to\R^d$ be a  piecewise regular path and let $k\in\N$. If $s,t\in [0,1]$, then
\begin{equation}\label{eq:holder}
\left| \sigma^{(k)}(X_{|[s,t]})\right| \apprle |t-s|^k.
\end{equation}
	\begin{proof}
		Let $S=\sigma^{(k)}(X_{|[s,t]})$. Since $k$ is fixed, in order to conclude it is enough to bound every entry of $S$. By definition
		\begin{align}
		|S_{i_1\ldots i_k}|&=\left|  \int_{s}^{t}\int_{s}^{t_k}\dots\int_{s}^{t_3}\int_{s}^{t_2}\dot{X}_{i_1}(t_1)\cdot\ldots\cdot\dot{X}_{i_k}(t_k)dt_1\dots dt_k\right| \nonumber\\
		&\le\int_{s}^{t}\int_{s}^{t_k}\dots\int_{s}^{t_3}\int_{s}^{t_2}\left| \dot{X}_{i_1}(t_1)\cdot\ldots\cdot\dot{X}_{i_k}(t_k)\right| dt_1\dots dt_k\nonumber\\
		&\le\sup_{t_1\in[0,1]}| \dot{X}_{i_1}(t_1)| \cdot\ldots\cdot\sup_{t_k\in[0,1]}| \dot{X}_{i_k}(t_k)|\int_{s}^{t}\int_{s}^{t_k}\dots\int_{s}^{t_3}\int_{s}^{t_2} dt_1\dots dt_k\nonumber\\
		&=\sup_{t_1\in[0,1]}| \dot{X}_{i_1}(t_1)| \cdot\ldots\cdot\sup_{t_k\in[0,1]}| \dot{X}_{i_k}(t_k)|\cdot\frac{|t-s|^k}{k!}.\nonumber
		\end{align}
		Since $X$ is  piecewise regular, all the suprema are finite.
	\end{proof}
\end{lem}

So we see that a  piecewise regular path $X:[0,1]\to\R^d$ induces a path $\sigma(X_{|[0,\textendash]}):[0,1]\to\G(\R^d)$ satisfying inequality (\ref{eq:holder}). If we want a rough path to be a generalization of a  piecewise regular path, we can define it in a similar flavor, also allowing different exponents. Recall that $p_k:T((\R^d))\to(\R^d)^{\otimes k}$ is the projection.

\begin{dfn}\label{dfn:roughpaths}
	A \textit{rough path} of order $m$ is a path ${\bf X}:[0,1]\to\G^m(\R^d)$ such that $|p_k({\bf X}(s)^{-1}\otimes {\bf X}(t))|\apprle |t-s|^\frac{k}{m}$
	for every $k\in\{1,\dots,m\}$ and every $s,t\in [0,1]$. The inverse is taken in the group $\G^m(\R^d)$.
\end{dfn}

In the literature, there are other classes of paths that are called rough, see for instance \cite[Definition 9.15]{Multidim} and \cite[Section 2.6]{Friz}. However, in this paper we only consider rough paths as in Definition \ref{dfn:roughpaths}, sometimes called \textit{weakly geometric H\"older} $m$-rough paths. Moreover, following \cite[Section 5.4]{AFS}, we will focus on a special subclass of rough paths of order $m$, indexed by elements of $\Lie^m(\R^d)$.

\begin{dfn}
	For $L\in\Lie^m(\R^d)$, consider the path ${\bf X}_L:[0,1]\to\G^m(\R^d)$ sending $t$ to $\exp^{(m)}(tL)$. %Define $${\bf X}_L=(p_1(g_L),p_2(g_L)).$$
	By \cite[Exercise 9.17]{Multidim}, this is indeed an order $m$ rough path. In analogy with the piecewise regular case, we define the \textit{signature} of ${\bf X}_L$ to be its endpoint $\sigma({\bf X}_L)=\exp(L)\in\G(\R^d)$.
\end{dfn}

We want to study the set parameterizing the $k$-th signatures of ${\bf X}_L$, when $L$ ranges over $\Lie^m(\R^d)$. Such set is the image of $p_k\circ\exp$. A priori, this is just a semialgebraic subset of $(\R^d)^{\otimes k}$, that is, it is described by a finite number of polynomial equations and inequalities. Semialgebraic sets are usually hard to handle. In order to make our analysis simpler, we will follow a common approach in applied algebraic geometry and take the Zariski closure of this set, which means that we only look at the equations without considering the inequalities. Furthermore, from a geometric viewpoint it is convenient to work over an algebraically closed field, so we look at the variety that these equations define in $(\C^d)^{\otimes k}$, instead of $(\R^d)^{\otimes k}$. Finally, we want to work up to scalar multiples, so we pass to the projectivization and hence we deal with a projective variety.

\begin{dfn} Let $f_{d,k,m}$ be the complexification of the map
	\[p_k\circ\exp:\Lie^m(\R^d)\rightarrow(\R^{d})^{\otimes k}.\]
	The \textit{Rough Veronese variety} $\RR_{d,k,m}\subset\p^{d^k-1}$ is the projectivization of the closure of the image of $f_{d,k,m}$. %:\Lie^m(\R^d)\xrightarrow{\exp}\G(\R^d)\xrightarrow{p_k}(\R^{d})^{\otimes k}\rightarrow (\R^{d})^{\otimes k}\otimes\C=(\C^{d})^{\otimes k}.\]	
\end{dfn}

First of all observe that, being the image of a map, $\RR_{d,k,m}$ is irreducible. There are several reasons to compare $\RR_{d,k,m}$ to a Veronese variety. Since $\Lie^1(\R^d)=\R^d$, an element $L\in\Lie^1(\R^d)$ is just a vector and therefore, up to a multiplicative constant, $p_k(\exp (L))=L^{\otimes k}$ can be viewed as the Veronese embedding of $\R^d$ into $\Sym^k\R^d\subset (\R^{d})^{\otimes k}$. In other words, $\RR_{d,k,1}=\V_{d-1,k}$. Moreover, $\V_{d,k}$ is toric, parametrized by all degree $k$ monomials in $d+1$ variables. We will see that $\RR_{d,k,m}$ is toric as well, and it is defined by monomials of weighted degree $k$, for a suitable choice of weights imposed by the structure of $\Lie^m(\R^d)$. Unlike $\V_{d,k}$, however, in general $\RR_{d,k,m}$ fails to be smooth.

The inclusions $\Lie^i(\R^d)\subset\Lie^{i+1}(\R^d)$ show that the Rough Veronese varieties are nested%, as $m$ increases
. On the other hand, this chain stabilizes. Indeed, when we apply $p_k$ and project onto the order $k$ summand, we do not see anything of order greater than $k$. So $\RR_{d,k,k+i}=\RR_{d,k,k}$ for every $i\in\N$. Being the image of $\G^k(\R^d)$ under $p_k$, $\RR_{d,k,k}$ contains all the $k$-th signatures of piecewise regular paths. For this reason, in \cite[Section 4.3]{AFS} it is called \textit{universal variety} and denoted by $\U_{d,k}$. Summarizing, there is a chain of strict inclusions
\[\V_{d-1,k}=\RR_{d,k,1}\subset\RR_{d,k,2}\subset\ldots\subset\RR_{d,k,k}=\RR_{d,k,k+1}=\ldots=\U_{d,k}.\]

It is not restrictive to assume $d\ge 2$ and $m\le k$. The first thing we want to do is to determine $\dim\RR_{d,k,m}$.

\begin{pro}
	If $m\le k$, then the general point of $\RR_{d,k,m}$ has $k$ preimages under $f_{d,k,m}$. In particular, $\dim\RR_{d,k,m}=\dim(\p(\Lie^m(\R^d)))$.
\end{pro}

This was already noted in \cite[Remark 6.5]{AFS}. We now want to understand the geometry of $\RR_{d,k,m}$, and we start by looking at the simplest example.

\begin{es}\label{ex:R222}
	Consider $d=k=m=2$. We want to write down $f_{2,2,2}$. The Lyndon words of length at most 2 in the alphabet $\{{\bf 1},{\bf 2}\}$ are {\bf 1}, {\bf 2} and {\bf 12}, hence $\dim\Lie^2(\R^2)=3$. An element of $\Lie^2(\R^2)$ can therefore be written as $x_1{\bf 1}+x_2{\bf 2}+a({\bf 12}-{\bf 21})$. If we look at it in tensor terms, we see a vector $x=(x_1,x_2)$ as order 1 summand, and a $2\times 2$ matrix $A$ as order 2 summand. Being a multiple of $({\bf 12}-{\bf 21})$, $A$ is skew-symmetric. Then
	\begin{align}
	%\pi_2(\exp(x_1{\bf 1}+x_2{\bf 2}+a({\bf 12}-{\bf 21})))&=
	p_2(\exp(x+A)) %\nonumber\\
	&=p_2\left(1+x+A+\frac{(x+A)^2}{2}+\ldots\right)\nonumber\\
	&=p_2\left(1+x+A+\frac{x^2+xA+Ax+A^2}{2}+\ldots\right)\nonumber\\
	&=A+\frac{x^2}{2}=\left(\begin{array}{cc}
	0 & a\\
	-a & 0\end{array}\right)+\frac{1}{2}\left(\begin{array}{cc}
	x_1^2 & x_1x_2\\
	x_1x_2 & x_2^2\end{array}\right).\nonumber
	\end{align}
	Then we can write $f_{2,2,2}:\C^3\rightarrow\C^{2\times 2}$ as
	\[(x_1,x_2,a)\mapsto\left(\begin{array}{cc}
	\frac{x_1^2}{2} & \frac{x_1x_2}{2}+a\\
	\frac{x_1x_2}{2}-a & \frac{x_2^2}{2}\end{array}\right).\]
	Up to a linear change of coordinates, it becomes
	\[(x_1,x_2,a)\mapsto\left(\begin{array}{cc}
	x_1^2 & x_1x_2\\
	a & x_2^2\end{array}\right).\]
	We can make several important remarks. First of all, the map is defined by monomials, so $\RR_{2,2,2}$ is a toric variety. It is a cone over the Veronese variety $\V_{1,2}$ and it spans the whole $\p^3$. Its ideal is generated by one quadratic polynomial. Finally, it can be seen as the embedding
	\[\p(1,1,2)\hookrightarrow\p^3
	\]
	of a weighted projective plane, defined by all monomials of (weighted) degree 2.
\end{es}
Our main goal is to generalize these remarks to all values of $d,k,m$. %The first task will be to prove that it is always possible to find a linear change of coordinates such that $\RR_{d,k,m}$ is toric.

\section*{$\RR_{d,k,m}$ as a toric variety}
Roughly speaking, a variety is toric if it is the image of a monomial map. A toric variety not only has nice properties - for instance, it is irreducible, rational and its ideal is generated by binomials - but it can be associated to a polytope that completely encodes its geometry. This makes toric varieties accessible from a theoretical, combinatorial and computational viewpoint. A good reference on toric varieties is \cite{cox}. This section is devoted to the proof that the Rough Veronese is indeed toric and to provide an explicit way to make computations on it.

If we write an element $L\in\Lie^m(\R^d)$ as the sum $L=L_1+\ldots+L_m$ of terms of order $1,\ldots,m$, then $p_k(\exp(L))\in(\R^d)^{\otimes k}$ is a linear combination of all possible ways to get an order $k$ tensor by multiplying $L_1,\ldots,L_m$. If we consider coordinates $T_{i_1\ldots i_k}$ on $(\C^d)^{\otimes k}$, we can rephrase this observation by saying that every coordinate $T_{i_1\ldots i_k}$ of $f_{d,k,m}$ is a linear combination of weighted degree $k$ monomials. The weight of a variable corresponding to a length $i$ Lyndon word is $i$. We can define a map
\[\Lie^m(\R^d)\to(\R^d)^{\otimes k}\]
by using all such weighted monomials, and consider its complexification $g_{d,k,m}$. By our observation, there is a linear change of coordinates sending the image of $g_{d,k,m}$ to the image of $f_{d,k,m}$. We want to prove that such change of coordinates is invertible, that is, every weighted monomial can be obtained as a linear combination of the entries of $f_{d,k,m}$, and we also want to do so as explicitly as possible. This section is based on %\cite[Section 3]{Joscha's notes} and
 \cite[Section IV]{exponential}. Let us start with two definitions.

%Consider the basis $\{P_{w}\mid w\in W_{d,m}\}$ of $\Lie^m(\R^d)$ indexed by the Lyndon words. Denote by $\iota:\Lie^m(\R^d)\to T(\R^d)$ the inclusion. We want to find elements $S_{w}\in T(\R^d)$ such that $\{\iota^*(S_{w})\mid w\in W_{d,m}\}$ is the dual basis of $\{P_{w}\mid w\in W_{d,m}\}$.

\begin{dfn}\label{def:S} Let $v$ be a word. If $v={\bf e}$ is the empty word, set $S_{{\bf e}}={\bf e}$. Otherwise, we define $S_v$ in the following recursive way.
	\begin{enumerate}
		\item If $v$ is Lyndon, write $v=l\cdot w$, where $l$ is a letter, and define $S_{v}=l\cdot S_w$.
		\item Otherwise, write $v=w_1^{ i_1}\cdot\ldots\cdot w_k^{i_k}$ as concatenation of decreasing Lyndon words. This can be done uniquely by \cite[Section 7.4]{bible}. Define
		\[S_v=\frac{S_{w_1}^{\shuffle i_1}\shuffle\ldots\shuffle S_{w_k}^{\shuffle i_k}}{i_1!\ldots i_k!}.\]
	\end{enumerate}
\end{dfn}

%By \cite[Theorem 5.3]{bible}, $\{S_{w}\mid w\in W_{d,m}\}$ have the required property.
The next ingredient we need is the following.

\begin{dfn}\label{def:pi}
	Define a linear map $\psi:T(\R^d)\to T(\R^d)$ that acts on a word $v$ by
	\[v\mapsto\sum_{n=1}^{\infty}\frac{(-1)^{n+1}}{n}\sum_{\substack{u_1,\dots,u_n \\ \mbox{\begin{tiny}
				nonempty
				\end{tiny}}
			\\ \mbox{\begin{tiny}
				words
				\end{tiny}}}} \langle v,u_1\cdot\ldots\cdot u_n\rangle u_1\shuffle\ldots\shuffle u_n.\]
\end{dfn}

In \cite{exponential}, the map $\psi$ is called $\pi'_1$%, while it appears in \cite{Joscha's notes} as $\pi^\top_1$
. Observe that only finitely many terms of the sum are non-zero. Indeed, an element in $T(\R^d)$ is the linear combination of finitely many words, and if $v$ is a word then there are only finitely many ways to write it as a concatenation $u_1\cdot\ldots\cdot u_n$. Moreover, if $v$ is a word of length $l$ and $v=u_1\cdot\ldots\cdot u_n$, then $u_1\shuffle\ldots\shuffle u_n$ is a sum of length $l$ words. This means that $\psi$ preserves the grading of $T(\R^d)$.

The following result about the exponential map will be of great help. Recall that $W_{d,m}$ is the set of Lyndon words of length at most $m$ in the alphabet with $d$ letters.

\begin{lem}
	\label{lem:explicit change of coord}
	Let $\{P_{w}\mid w\in W_{d,m}\}$ be the basis of $\Lie^m(\R^d)$ indexed by the Lyndon words and let $w_1,\dots,w_r\in W_{d,m}$. Then
	\[\left\langle\exp\left( \sum_{w\in W_{d,m}}\alpha_wP_{w}\right) ,\psi(S_{w_1})\shuffle\ldots\shuffle \psi(S_{w_r})\right\rangle=\alpha_{w_1}\cdot\ldots\cdot \alpha_{w_r}.\]
	\begin{proof}
		By \cite[Theorem 1]{exponential}, for every $i\in\{1,\dots,r\}$ we have
		\[\left\langle\exp\left( \sum_{w\in W_{d,m}}\alpha_wP_{w}\right) ,\psi(S_{w_i})\right\rangle=\alpha_{w_i}.\]
	Now the statement follows by Lemma \ref{lem:shuffle identity}.
	\end{proof}
\end{lem}

Clearly every entry of $f_{d,k,m}$ is a linear combination of monomials. Lemma \ref{lem:explicit change of coord} shows that the converse holds. Every monomial can be obtained as linear combination of entries of $f_{d,k,m}$. In other words, we can use the map $\psi$ to build the linear forms we needed to pass from $f_{d,k,m}$ to $g_{d,k,m}$, allowing us to identify them. Let us summarize the conclusion.

\begin{pro}\label{pro:toric} %Set $\mu=\mu_{1,d}+\ldots+\mu_{m,d}$ and let $w_1,\ldots,w_\mu$ be the Lyndon words of length at most $m$
%Let $W_m$ be the set of Lyndon words of length at most $m$. 
For every $w\in W_{d,m}$, define a variable $x_w$ and assign it the weight $|w|$. Then, up to a linear change of coordinates, $f_{d,k,m}$ is defined by all monomials of weighed degree $k$ in the variables $\{x_w\mid w\in W_{d,m}\}$. More explicitly, if we set $$J=\{(w_1,\dots, w_r)\mid r\in\N, w_i\in W_{d,m}\mbox{ and } |w_1|+\ldots+|w_r|=k\},$$ then  %$f_{d,k,m}:\Lie^m(\R^d)\to\p^{d^k-1}$ can be written as
$\RR_{d,k,m}$ is isomorphic to the image of
\[[x_w]_{w\in W_{d,m}}\mapsto [x_{w_1}\cdot\ldots\cdot x_{w_r}]_{(w_1,\dots, w_r)\in J}.	\] %	followed by the inclusion in $\p^{d^k-1}$ by adding as many zeroes as necessary.
	In particular, it is a toric variety.
\end{pro}

\begin{rmk}
	We could phrase Proposition \ref{pro:toric} in a slightly different way by saying that, up to a linear change of coordinates, $p_k\circ\exp$ is a monomial map. More precisely, for every $d,k,m\in\N$ there exists a linear change of coordinates $\lambda:(\R^d)^{\otimes k}\to (\R^d)^{\otimes k}$ such that the map
	\[\lambda\circ p_k\circ \exp:\Lie^m(\R^d)\to (\R^d)^{\otimes k}\]
	is defined by all monomials of weighted degree $k$.
\end{rmk}

Notice that the image of the map defined in Proposition \ref{pro:toric} lies in a possibly proper linear subspace of $\p^{d^k-1}$. We will discuss the linear span of $\RR_{d,k,m}$ in Lemma \ref{pro:countmonomials}.

For the purpose of practical applications, we want to explicitly describe a linear change of coordinates in $\p^{d^k-1}$ that makes $\RR_{d,k,m}$ toric. If $m<k$, then every element of $L\in\Lie^{m+1}(\R^d)$ can be written uniquely as $L=(L_1,L_2)$, where $L_1\in\Lie^m(\R^d)$. In this case, $f_{d,k,m}(L_1)=f_{d,k,m+1}(L_1,0)$. This means that every change of coordinates that makes $f_{d,k,m+1}$ a monomial map also makes $f_{d,k,m}$ monomial. Therefore we can assume $m=k$, because the change of coordinates that will make $\RR_{d,k,k}$ toric will also work on $\RR_{d,k,m}$ for every $m\le k$.

In the notation of Proposition \ref{pro:toric}, % again we can define$$J=\{(w_1,\dots, w_r)\mid r\in\N, w_i\in W_{d,m}\mbox{ and } |w_1|+\ldots+|w_r|=k\}.$$ Then 
the change of coordinates in $\p^{d^k-1}$ is
\begin{align}\label{eq:changeofcoordinates}
T\mapsto \left[ \langle T,\psi(S_{w_1})\shuffle\ldots\shuffle \psi(S_{w_r})
\rangle\right]_{(w_1,\dots, w_r)\in J}.
\end{align}
This is indeed well defined because, as we will show in Proposition \ref{pro:nondegen}, $J$ has exactly $d^k$ elements when $m=k$.

\begin{es}
	\label{ex:change of coordinates U24} Let us consider $\U_{2,4}=\RR_{2,4,4}$ and let $T_{ijkl}$ be the coordinates of $(\R^2)^{\otimes 4}$. The set of Lyndon words of length at most 4 in the alphabet $\{{\bf 1},{\bf 2}\}$ is
	\[W_{2,4}=\{{\bf 1},{\bf 2},{\bf 12},{\bf 112},{\bf 122},{\bf 1112},{\bf 1122},{\bf 1222}\}.\]
	By using Definition \ref{def:S}, in this case we get $S_w=w$ for every $w\in W_{2,4}$.
	The first 5 entries of the change of variables (\ref{eq:changeofcoordinates}) are then
	\begin{align}
	\langle T,\psi(S_{{\bf i}})\shuffle\psi(S_{{\bf j}})\shuffle \psi(S_{{\bf k}})\shuffle \psi(S_{{\bf l}})
	\rangle&=\langle T,\psi({\bf i})\shuffle\psi({\bf j})\shuffle \psi({\bf k})\shuffle \psi({\bf l})
	\rangle\nonumber\\
	&=\langle T,{\bf i}\shuffle{\bf j}\shuffle {\bf k}\shuffle {\bf l}
	\rangle\nonumber\\
	&=\sum_{\sigma}T_{\sigma(i)\sigma(j)\sigma(k)\sigma(l)},\nonumber
	\end{align}
	where $1\le i\le j\le k\le l\le 2$ and $\sigma$ ranges among the permutations of $\{i,j,k,l\}$. The next 3 entries are of the form
	\begin{align}
	\langle T,\psi(S_{{\bf 12}})\shuffle\psi(S_{{\bf 1}})\shuffle \psi(S_{{\bf 1}})\rangle&=\langle T,\psi({\bf 12})\shuffle\psi({\bf 1})\shuffle \psi({\bf 1})\rangle\nonumber\\
	&=\left\langle T,\frac{1}{2}({\bf 12}-{\bf 21})\shuffle{\bf 1}\shuffle {\bf 1}\right\rangle \nonumber\\
	&=\langle T, 3\cdot {\bf 1112}+ {\bf 1121}- {\bf 1211}-3\cdot {\bf 2111}\rangle \nonumber\\
	&=3T_{1112}+T_{1121}-T_{1211}-3T_{2111}.\nonumber
	\end{align}
We continue in this way until we reach length 4 Lyndon words
\begin{align}
\langle T,\psi(S_{{\bf 1112}})\rangle=\langle T,\psi({\bf 1112})\rangle=\left\langle T,\frac{1}{6}({\bf 1211}-{\bf 1121})
\right\rangle=\frac{1}{6}(T_{1211}-T_{1121})\nonumber
\end{align}
and, in the same way,
\begin{align}
\langle T,\psi(S_{{\bf 1122}})\rangle&=\frac{1}{6}(T_{1122}-T_{1212}+T_{2121}-T_{2211}),\nonumber\\
\langle T,\psi(S_{{\bf 1222}})\rangle&=\frac16(T_{2212}-T_{2122}).\nonumber
\end{align}
Compare our computations to the ones in \cite[Section IV]{exponential}. If we denote by $Q_{ijkl}$ the new coordinates on $(\R^2)^{\otimes 4}$, then our change of coordinates in $\p^{15}$ is
\[\left[ \begin{matrix}
	Q_{1111}\\
	Q_{1112}\\
	Q_{1121}\\
	Q_{1122}\\
	Q_{1211}\\
	Q_{1212}\\
	Q_{1221}\\
	Q_{1222}\\
	Q_{2111}\\
	Q_{2112}\\
	Q_{2121}\\
	Q_{2122}\\
	Q_{2211}\\
	Q_{2212}\\
	Q_{2221}\\
	Q_{2222}\\
\end{matrix}\right]=\left[ \begin{matrix}
24T_{1111}\\
6(T_{1112}+T_{1121}+T_{1211}+T_{2111})\\
4(T_{1122}+T_{1212}+T_{1221}+T_{2211}+T_{2112}+T_{2121})\\
6(T_{1222}+T_{2122}+T_{2212}+T_{2221})\\
24T_{2222}\\
3T_{1112}+T_{1121}-T_{1211}-3T_{2111}\\
2T_{1122}+T_{1212}-T_{2121}-2T_{2211}\\
3T_{1222}+T_{2122}-T_{2212}-3T_{2221}\\
\frac12(T_{1112}-T_{1121}-T_{1211}+T_{2111})\\
\frac13(T_{1122}+T_{2112}+T_{2211})-\frac16(T_{1212}+T_{2121}+4T_{1221})\\
\frac13(T_{1122}+T_{1221}+T_{2211})-\frac16(T_{1212}+T_{2121}+4T_{2112})\\
\frac12(T_{1222}-T_{2122}-T_{2212}+T_{2221})\\
T_{1122}-T_{1221}-T_{2112}+T_{2211}\\
\frac{1}{6}(T_{1211}-T_{1121})\\
\frac{1}{6}(T_{1122}-T_{1212}+T_{2121}-T_{2211})\\
\frac16(T_{2212}-T_{2122})
\end{matrix}\right]. \]
\end{es}

In Example \ref{ex:change of coordinates U24}, every $S_w$ constructed via Definition \ref{def:S} coincides with the corresponding word $w$. However, this is not true in general, as the case $d=k=3$ shows.

\begin{es}\label{ex:change of coordinates U33} Let us consider $\U_{3,3}=\RR_{3,3,3}$% and let $T_{ijk}$ be the coordinates of $(\R^3)^{\otimes 3}$
	. The set of Lyndon words of length at most 3 in the alphabet $\{{\bf 1},{\bf 2},{\bf 3}\}$ is
	\[W_{3,3}=\{{\bf 1},{\bf 2},{\bf 3},{\bf 12},{\bf 13},{\bf 23},{\bf 112},{\bf 113},{\bf 122},{\bf 133},{\bf 223},{\bf 233},{\bf 123},{\bf 132}\}.
	\]
Here we observe that $S_{{\bf 132}}={\bf 123}+{\bf 132}$. On the other hand, $S_w=w$ for every $w\in W_{3,3}\setminus\{{\bf 132}\}$. By carrying out the computations as in Example \ref{ex:change of coordinates U24}, we find that the first 10 entries of the change of variables (\ref{eq:changeofcoordinates}) are
 \begin{align}
\langle T,\psi(S_{{\bf i}})\shuffle\psi(S_{{\bf j}})\shuffle \psi(S_{{\bf k}})
\rangle%&=\langle T,\psi({\bf i})\shuffle\psi({\bf j})\shuffle \psi({\bf k}) \rangle\nonumber\\
% &=\langle T,{\bf i}\shuffle{\bf j}\shuffle {\bf k} \rangle\nonumber\\
% &=\langle T,{\bf kij}+{\bf ikj}+{\bf ijk}+{\bf kji}+{\bf jki}+{\bf jik} \rangle\nonumber\\
 &=T_{kij}+T_{ikj}+T_{ijk}+T_{kji}+T_{jki}+T_{jik}\nonumber
 \end{align}
 for every $1\le i\le j\le k\le 3$. Then we have 9 more entries of the form
\begin{align}
\langle T,\psi(S_{{\bf i}})\shuffle\psi(S_{{\bf jk}})\rangle%&=\langle T,\psi({\bf i})\shuffle\psi({\bf jk})\rangle\nonumber\\
% &=\langle T,{\bf i}\shuffle\frac{1}{2}({\bf jk}-{\bf kj}) \rangle\nonumber\\
% &=\langle T,\frac{1}{2}({\bf ijk}+{\bf jik}+{\bf jki}-{\bf ikj}-{\bf kij}-{\bf kji})\rangle\nonumber\\
 &=\frac{1}{2}(T_{ijk}+T_{jik}+T_{jki}-T_{ikj}-T_{kij}-T_{kji}),\nonumber
\end{align}for every $i\in\{1,2,3\}$ and $1\le j<k\le 3$. Finally, we tackle length 3 Lyndon words. We have
\begin{align}
\langle T,\psi(S_{{\bf 132}})\rangle%&=\langle T,\psi({\bf 123}+{\bf 132})\rangle\nonumber\\
%&=\langle T,\psi({\bf 123})+\psi({\bf 132})\rangle\nonumber\\
%&=\langle T,\frac{1}{6}({\bf 123}+{\bf 132}-2\cdot{\bf 213}+{\bf 231}-2\cdot{\bf 312}+{\bf 321})\rangle\nonumber\\
&=\frac{1}{6}(T_{123}+T_{132}+T_{231}+T_{321})-\frac{1}{3}(T_{213}+ T_{312}),\nonumber
\end{align}while the last 7 entries are
\begin{align}
\langle T,\psi(S_{{\bf ijk}})\rangle%&=\langle T,\psi({\bf ijk})\rangle\nonumber\\
%&=\langle T,\frac{1}{6}(2\cdot{\bf ijk}+2\cdot{\bf kji}-{\bf kij}-{\bf ikj}-{\bf jki}-{\bf jik})\rangle\nonumber\\
&=\frac{1}{3}(T_{ijk}+T_{kji})-\frac{1}{6}(T_{kij}+T_{ikj}+T_{jki}+T_{jik}),\nonumber
\end{align}for ${\bf ijk}\in\{{\bf 112},{\bf 113},{\bf 122},{\bf 133},{\bf 223},{\bf 233},{\bf 123}\}$.
It is interesting to point out that some of these linear forms are used in \cite[Section 5]{annaemax} in order to recover the path of a given signature.\end{es}

\begin{rmk} The universal variety contains many interesting subvarieties besides $\RR_{d,k,m}$. Examples include the signature variety $\LL_{d,k,m}$ of piecewise linear paths with $m$ steps and its subvariety $\AAA_{\nu,k}$ of axis-parallel paths, both studied in \cite{AFS}. The change of coordinates given by Lemma \ref{lem:explicit change of coord} proves that $\RR_{d,k,m}$ and therefore $\U_{d,k}$ are toric, but it does not necessarily work as nicely with other subvarieties. For instance, a software computation shows that, after our change of coordinates, the ideals of both $\LL_{2,3,2}$ and $\AAA_{1212,3}$ are not generated by binomials. The question of whether $\AAA_{\nu,k}$ is toric is addressed in \cite[Section 3.1]{CGM}, but the answer is not known in full generality.
\end{rmk}

\section*{Further geometric properties}

In this section we will generalize the remarks we made for Example \ref{ex:R222}. While $\RR_{2,2,2}\subset\p^3$ is nondegenerate, for other values of $d,k,m$ the Rough Veronese may be contained in a smaller linear subspace, as we already observed in Proposition \ref{pro:toric}.

\begin{pro}\label{pro:countmonomials}
The affine dimension of the linear span of $\RR_{d,k,m}\subset\p^{d^k-1}$ is
	\[\sum_{\lambda\vdash k,\lambda_1\le m}\left( \prod_{i=1}^{m}\binom{\mu_{i,d}+\sharp\{j\mid\lambda_j=i\}-1}{\mu_{i,d}-1}\right) .\]
	This number equals the coefficient of $t^k$ of the expansion of the generating function
	\[\prod_{i=1}^{m}\frac{1}{(1-t^i)^{\mu_{i,d}}}.\]
	\begin{proof} As in Proposition \ref{pro:toric}, we define a variable $x_w$ for every $w\in W_{d,m}$ and we assign it the weight $|w|$.
		Given a weighted degree $k$ monomial $\f\in\C[x_w\mid w\in W_{d,m}]$, we can write it in reverse lexicographic order as a string of possibly repeated variables. Define a partition $\lambda$ of $k$ by
		\[\lambda_j=i\Leftrightarrow\mbox{ the $j$-th entry of the string is } x_w \mbox{ for some } |w|=i.\]
		Now $\f$ is the product of monomials $\f=\f_1\cdot\ldots\cdot \f_m$, where $\f_i$ is a monomial in $\C[x_w\mid |w|=i]$.
		The degree of $\f_i$ is the number of times $i$ appears as an entry of $\lambda$. Hence, for every $\f_i$ there are $\binom{\mu_{i,d}+\sharp\{j\mid\lambda_j=i\}-1}{\mu_{i,d}-1}$ choices. For the generating function, see \cite[Section 1.8]{combinat}.
	\end{proof}
\end{pro}

Nonetheless, the universal variety $\U_{d,k}=\RR_{d,k,k}$ is indeed nondegenerate for every $d$ and $k$.

\begin{pro}\label{pro:nondegen} Under our assumption $k\ge m$, the rough Veronese $\RR_{d,k,m}\subset\p^{d^k-1}$ is nondegenerate if and only if $m=k$.
	\begin{proof}
		By Proposition \ref{pro:countmonomials}, the dimension of the linear span of $\RR_{d,k,m}$ is strictly increasing with respect to $m$, hence % than that of $\RR_{d,k,k}$. Because of the chain of strict inclusions $\RR_{d,k,m}\subsetneq\RR_{d,k,m+1}$,
		 it is enough to show that $\RR_{d,k,k}$ is nondegenerate.
		
		Thanks to Proposition \ref{pro:toric}, we only have to prove that there are $d^k$ distinct monomials of weighted degree $k$. Proposition \ref{pro:countmonomials} may be difficult to apply, so instead we want to define a bijection between these monomials and the set of length $k$ words in the alphabet with $d$ letters. Since our variables are indexed by Lyndon words, we can think of a monomial as a product of Lyndon words such that the sum of the lengths is $k$. Observe that the monomial remains the same if we permute the variables, i.e.~the Lyndon words. In order to avoid redundancy, we can fix an order. Now we only have to prove that, after fixing an order among Lyndon words, every length $k$ word in the alphabet with $d$ letters can be written uniquely as an ordered product of Lyndon words whose lengths sum to $k$. This is a well-known fact, and a proof can be found for instance in \cite[Section 7.4]{bible}.
	\end{proof}
\end{pro}

The next result shows another feature of rough paths. Not only does the universal variety $\U_{d,k}$ coincide with the last Rough Veronese $\RR_{d,k,k}$, but its structure is already determined by the second to last one, $\RR_{d,k,k-1}$. 

\begin{pro}\label{pro:universal}
	The universal variety $\U_{d,k}=\RR_{d,k,k}$ is a cone over $\RR_{d,k,k-1}$ with vertex $\p^{\mu_{k,d}-1}$. The preimage of the vertex is defined by the vanishing of the first $k-1$ entries.
	\begin{proof}
	Let $V=\R^{\mu_{d,k}}$ be the vector subspace of $\Lie^k(\R^d)$ defined by the vanishing of the first $k-1$ entries. Thanks to Proposition \ref{pro:toric}, we know that $\RR_{d,k,k}$ is defined by the degree $k$ weighted monomials. Exactly $\mu_{k,d}$ of them are $\{x_w\mid w\in W_{d,m}\mbox{ and } |w|=k\}$. In other words, the map $p_k\circ\exp$ restricts to the identity on $V$. It follows that $$p_k(\exp(\Lie^k(\R^d)))=p_k(\exp(\Lie^{k-1}(\R^d)))\times V$$ is a cylinder. When we pass to the projectivization, we get a cone with base $\RR_{d,k,k-1}$ and vertex $\p(p_k(\exp(V))\otimes\C)=\p(V\otimes\C)$.
	\end{proof}
\end{pro}

\begin{rmk}	It is interesting to try to classify the paths whose signatures lie in the vertex of the universal variety. Equivalently, we wonder what it means for a rough path to have zeroes in the first $k-1$ entries. It is straightforward to check that for a smooth path $X$, the first signature
\[\sigma^{(1)}(X)=X(1)-X(0)\]
is just the vector joining the endpoint of $X$ to its starting point. The second signature also has a geometric interpretation: if we take the projection of $X$ onto the $i,j$ plane, the signed area of the region bounded by $X$ and the segment between $X(0)$ and $X(1)$ is  $\frac{1}{2}(\sigma(X)_{ij}-\sigma(X)_{ji})$. For instance, if $X$ is a loop then $\sigma^{(1)}(X)=0$ but $\sigma^{(2)}(X)$ may be nonzero. However, if we allow not only smooth paths, but also rough paths, we find more interesting examples. There are rough paths ${\bf X}$ such that $p_1({\bf X}(t))=p_1({\bf X}(0))$ for every $t$, but with a nonzero second entry. If we try to imagine it, our intuition suggests a constant path that nonetheless describes an area. This is not possible for a piecewise regular path. One example of these \textit{pure area paths} is \cite[Exercise 2.17]{Friz}, where it is built as a limit of smaller and smaller loops.

In the same fashion, we can think of the vanishing of the first two signatures as a loop with some symmetry that makes the signed area zero, such as two circles in $\R^2$ meeting at a point. Then again, there are more exotic examples of rough paths with the same property. However, the geometric meaning of the third and higher order signatures is not yet clear, as it is an open research area.
\end{rmk}

We now turn our attention to the ideal of $\RR_{d,k,m}$. In \cite[Section 4]{annaemax}, it is proved that the ideal of $\RR_{d,3,3}$ is generated by the $2\times 2$ minors of a suitable Hankel matrix, in particular by quadrics. In many other examples this ideal is generated in degree 2. However, this is not true in general, not even if $k=m$. The following counterexample, suggested by Micha\l{}ek, answers a question posed in \cite[Section 4.3]{AFS}.

\begin{pro}\label{pro:not generated by quadrics} For $14\le m\le k=20$, the ideal of the Rough Veronese variety $\RR_{d,20,m}$ is not generated by quadrics.
	\begin{proof}
		Let $N$ be the number of weighted degree $k$ monomials. By Proposition \ref{pro:toric}, up to change of coordinates $\RR_{d,20,m}$ lies in $\p^{N-1}$. The coordinates $T$ of $\p^{N-1}$ are indexed by the $N$ weighted monomials. For every $i\in\{1,\dots,m\}$, let $x^{(i)}$ be one of the variables of weight $i$. Let $I$ be the ideal of $\RR_{d,20,m}$ and define \begin{align}
		t_1&=T_{x^{(1)}x^{(9)}x^{(10)}},& t_2 &= T_{x^{(5)}x^{(7)}x^{(8)}},\nonumber\\
		t_3&=T_{x^{(2)}x^{(4)}x^{(14)}}, & t_4&=T_{x^{(1)}x^{(5)}x^{(14)}},\nonumber\\
		t_5&=T_{x^{(4)}x^{(7)}x^{(9)}},& t_6&=T_{x^{(2)}x^{(8)}x^{(10)}}.\nonumber
		\end{align}
		Then $f=t_1t_2t_3-t_4t_5t_6$ is a degree 3 element of $I$. Let us show that $f$ is not generated by quadrics. Let $q_1,\dots,q_r$ be the degree two generators of $I$. Since $\RR_{d,20,m}$ is toric, they are binomials and so we can write $q_i=g_i-h_i$, where each $g_i$ and each $h_i$ is a degree 2 monomial. Assume by contradiction that $f$ can be algebraically generated by $q_1,\dots,q_r$. Then there exist degree 1 polynomials $l_1,\dots,l_r$ such that, up to order, $f$ is a sum of nonzero polynomials $f=l_1(g_1-h_1)+\ldots+l_t(g_r-h_r)$. It follows that each term in the sum is a multiple of some $g_i$ or some $h_i$. For instance, there is a degree 2 monomial dividing $t_1t_2t_3$. The only three possibilities are $t_1t_2$, $t_1t_3$ and $t_2t_3$. Suppose then that $g_i=t_1t_2$. Since $g_i-h_i\in I$, the product of the two monomials indexing the variables of $h_i$ equals $t_1t_2$. But it is easy to see that $1+9+10$ and $5+7+8$ are the only ways to obtain 20 by sums of non-repeated elements of $\{1,5,7,8,9,10\}$. Hence $g_i-h_i=0$, a contradiction.
	\end{proof}
\end{pro}

Among the features of $\RR_{2,2,2}$ we pointed out in Example \ref{ex:R222}, there is one we still have to check. We saw that $\RR_{2,2,2}$ is the embedding of $\p(1,1,2)$ given by its weighted quadrics. If we define the sequence of weights
\[s=(\underbrace{1,\dots,1}_{\mu_{1,d}\mbox{ {\Small times}}},\underbrace{2,\dots,2}_{\mu_{2,d}\mbox{ {\Small times}}},\dots,\underbrace{m,\dots,m}_{\mu_{m,d}\mbox{ {\Small times}}}),
\]
then by Proposition \ref{pro:toric} $f_{d,k,m}$ induces a rational map $\p(s)\dashrightarrow\p^{d^k-1}$, that we will still call in the same way by abuse of notation. However, $f_{d,k,m}$ does not need to be an embedding. Actually, it does not even need to be defined everywhere. For instance, $f_{2,3,2}:\p(1,1,2)\dashrightarrow\p^7$ is defined by
\[[x_1,x_2,a]\mapsto \left[\begin{array}{cc}
x_1^3 & x_1^2x_2\\
x_1x_2^2 & x_2^3\\\hline
x_1a & x_2a\\
0 & 0\end{array}\right],\]
so $[0,0,1]$ is a base point. This is a general behavior.

\begin{pro}
	The map $f_{d,k,m}:\p(s)\dashrightarrow\p^{d^k-1}$ is base point free if and only if every entry of $s$ divides $k$.
	\begin{proof}
	Assume that $k$ is a multiple of every entry of $s$. Then for every variable $x_w$, there is a power of $x_w$ appearing among the monomials of weighted degree $k$. Therefore the only way for all weighted monomials to vanish is setting all variables to zero. This means there are no base points.
		
	On the other hand, assume that $k$ is not a multiple of one of the entries of $s$, say $i$, and consider a variable $x^{(i)}$ of weight $i$. Since $i\nmid k$, no power of $x^{(i)}$ appears among the monomials defining $f_{d,k,m}$. Then $[0,\dots,0,1,0\dots,0]\in\p(s)$, with a 1 in the entry corresponding to $x^{(i)}$, is a base point for $f_{d,k,m}$.
	\end{proof}
\end{pro}

\section*{Table of invariants}\label{sec:table}
We collect some of the geometric invariants of $\RR_{d,k,m}$, obtained with the software Macaulay2. We compute the projective dimension of the linear span, the dimension of $\RR_{d,k,m}$ and its degree. Despite Proposition \ref{pro:not generated by quadrics}, all the ideals in the examples we present are generated by quadrics. In the last column \textquotedblleft gen\textquotedblright\ we record the number of generators.
\begin{center}
	\begin{tabular}{|c c c %|c
			|c c c|c|}
		\hline
		$d$	&	$k$	&	$m$	%& map
		 &	span   & $\dim$ & $\deg$	& gen \\ \hline %{|===|=|===|=|}
		$2$	&	$2$	&	$2$	%& $\R^2\times\R\to\R^4$
			  & 3     & 2      & 2    	& 1   \\
		$3$	&	$2$	&	$2$	%&	$\R^3\times\R^3\to\R^9$
		  & 8     & 5      & 4    	& 6   \\
		
		$4$	&	$2$	&	$2$	%&	$\R^4\times\R^6\to\R^{16}$
		  & 15   & 9     & 8    	& 20   \\
		
		$5$	&	$2$	&	$2$	%&	$\R^5\times\R^{10}\to\R^{25}$ 
		 & 24    & 14     & 16    	& 50   \\
		
		$6$	&	$2$	&	$2$	%&	$\R^6\times\R^{15}\to\R^{36}$
		  & 35    & 20     & 32    	& 105   \\
%		\hline
		$2$	&	$3$	&	$2$	%&	$\R^2\times\R\to\R^8$
		  & 5     & 2      & 4    	& 6   \\
		
		$3$	&	$3$	&	$2$	%&	$\R^3\times\R^3\to\R^{27}$ 
		 & 18    & 5    & 24    	& 81   \\
		
		$4$	&	$3$	&	$2$	%&  $\R^4\times\R^6\to\R^{64}$
		 &	43    & 9     & 200   	& 486   \\
		
		$5$	&	$3$	&	$2$	%&	$\R^5\times\R^{10}\to\R^{125}$
		  & 84    & 14     &   2221  	& 1920      \\
%		\hline
		$2$	&	$4$	&	$2$	%&	$\R^2\times\R\to\R^{16}$
		  & 8     & 2      & 8    	& 27   \\
		
		$3$	&	$4$	&	$2$	%&	$\R^3\times\R^3\to\R^{81}$
		  & 38   &   5    &  128   	& 528   \\
%		$4$	&	$4$	&	$2$	&	$\R^4\times\R^6\to\R^{256}$  &  115  &    9  &     	&   \\
%		\hline
		$2$	&	$5$	&	$2$	%&	$\R^2\times\R\to\R^{32}$
		  & 11   &    2 &   12  	&  43 \\
		
		$3$	&	$5$	&	$2$ %	& $\R^3\times\R^3\to\R^{243}$	
		  & 68   &   5  &   368  	&  1806\\
%		\hline
		$2$	&	$6$	&	$2$	%&	$\R^2\times\R\to\R^{64}$ 
		 &  15  &    2 &  18  	&  87 \\
		\hline %{|===|=|===|=|}
		$2$	&	$3$	&	$3$	%&	$\R^2\times\R\times\R^2\to\R^{8}$
		  & 7     &  4     &  4  	& 6   \\
		
		$3$	&	$3$	&	$3$	%&$\R^3\times\R^3\times\R^8\to\R^{27}$
			 & 26   &   13    &  24   	&   81   \\
		
		$4$	&	$3$	&	$3$	%& $\R^4\times\R^6\times\R^{20}\to\R^{64}$
			 &  63  &   29    &    200 	&   486   \\
	%	\hline
		$2$	&	$4$	&	$3$	%&	$\R^2\times\R\times\R^2\to\R^{16}$ 
		 &  12 &   4    &   12 	&  33  \\
		
		$3$	&	$4$	&	$3$	%&	$\R^3\times\R^3\times\R^8\to\R^{81}$  
		&  62  &   13   &  672  	& 954 \\
		
		$2$	&	$5$	&	$3$	&  19  &   4  &  28  	& 108 \\
		\hline %{|===|=|===|=|}
		$2$	&	$4$	&	$4$	%&	$\R^2\times\R\times\R^2\times\R^3\to\R^{16}$
		  & 15   &   7    &   12 	&   33 \\
%		$3$	&	$4$	&	$4$	&	  &    &      &   	& \\
%		\hline
		$2$	&	$5$	&	$4$ %	&	$\R^2\times\R\times\R^2\times\R^3\to\R^{32}$ 
		 &  25  &    7   &   40 	&  150\\
%		$3$	&	$5$	&	$4$	&$\R^3\times\R^3\times\R^8\times\R^{18}\to\R^{243}$	  &  194  &    31  &   	& \\
		\hline %{|===|=|===|=|}
		$2$	&	$5$	&	$5$	%&	$\R^2\times\R\times\R^2\times\R^3\times\R^6\to\R^{32}$ 
		 &   31 &   13    &  40  	&  150\\
%		\hline
		$2$	&	$6$	&	$5$	%&	$\R^2\times\R\times\R^2\times\R^3\times\R^6\to\R^{64}$  
		&   54&   13    &  336  	& 694 \\
		\hline
	\end{tabular}
\end{center}

\end{document}